\documentclass{amsart}
\usepackage{amssymb}
\usepackage{latexsym}
\usepackage{amsmath}
\usepackage{euscript}

\def\sD{{\mathfrak D}}      
\def\sG{{\mathfrak G}}   \def\sH{{\mathfrak H}}   
   \def\sK{{\mathfrak K}}   
      
      \def\sR{{\mathfrak R}}

      \def\dC{{\mathbb C}}

\def\bm\chi{\mbox{\boldmath$\chi$}}
\def\half{{\frac{1}{2}}}

\def\ker{{\rm ker\,}}
\def\ran{{\rm ran\,}}
\def\cran{{\rm \overline{ran}\,}}
\def\dom{{\rm dom\,}}
\def\mul{{\rm mul\,}}

\def\cdom{{\rm \overline{dom}\,}}

\def\okr{\stackrel{\scriptscriptstyle{\sf{def}}}{=}}

\let\xker=\ker \def\ker{{\xker\,}}
\def\span{{\rm span\,}}

\def\uphar{{\upharpoonright\,}}

\newtheorem{theorem}{Theorem}[section]
\newtheorem{proposition}[theorem]{Proposition}
\newtheorem{corollary}[theorem]{Corollary}
\newtheorem{lemma}[theorem]{Lemma}

\DeclareMathOperator{\hplus}{\, \widehat + \,}
\DeclareMathOperator{\hoplus}{\, \widehat \oplus \,}
\DeclareMathOperator{\hominus}{\, \widehat \ominus \,}

\theoremstyle{definition}
\newtheorem{example}[theorem]{Example}
\newtheorem{remark}[theorem]{Remark}

\numberwithin{equation}{section}
\begin{document}
\hfill\small{to appear in Acta Math. Hungarica July 2007, volume 116(1-2)}
\vspace{1cm}
\title[A canonical decomposition for linear relations]
{A canonical decomposition for linear operators and linear
relations}
\author{S.~Hassi}
\author{Z.~Sebesty\'en}
\author{H.S.V.~de~Snoo}
\author{F.H.~Szafraniec}
\address{Department of Mathematics and Statistics \\
University of Vaasa \\
P.O. Box 700, 65101 Vaasa \\
Finland}
\email{sha@uwasa.fi}
\address{Department of Applied Analysis \\
E\"otv\"os Lor\'and University\\
P\'azm\'any P\'eter s\'et\'any 1/C \\
1117 Budapest\\
Hungary}
\email{sebesty@cs.elte.hu}
\address{Department of Mathematics and Computing Science\\
University of Groningen \\
P.O. Box 800, 9700 AV Groningen \\
Nederland}
\email{desnoo@math.rug.nl}
\address{Institute of Mathematics \\
Uniwersytet Jagiello\'nski \\
ul. Reymonta 4 \\
30059 Krak\'ow \\
Poland}
\email{fhszafra@im.uj.edu.pl}
\thanks{The research was supported by bilateral agreements between
the E\"otv\"os Lor\'and University in Budapest, the University of
Groningen, and the Jagiellonian University in Krakow, and by the
Research Institute for Technology of the University of Vaasa. The
fourth author was also supported by the KBN grant 2 PO3A 037 024.}



\begin{abstract}
An arbitrary linear relation (multivalued operator) acting from one
Hilbert space to another Hilbert space is shown to be the sum of a
closable operator and a singular relation whose closure is the
Cartesian product of closed subspaces. This decomposition can be
seen as an analog of the Lebesgue decomposition of a measure into a
regular part and a singular part. The two parts of a relation are
characterized metrically and in terms of Stone's characteristic
projection onto the closure of the linear relation.
\end{abstract}

\maketitle

\section{Introduction}

Let $T$ be a linear operator from a Hilbert space $\sH$ to a Hilbert
space $\sK$. In general the closure $\overline{T}$ of the operator
$T$ (i.e., the closure of the graph of $T$ in the Cartesian product
$\sH \times \sK$) is not the graph of an operator anymore. In other
words, (the graph of) $\overline{T} $ has a nontrivial multivalued
part $\mul \overline{T}=\{\,k \in \sK :\, \{0,k\} \in
\overline{T}\,\}$. Relative to the closed linear subspace $\mul
\overline{T}$ of $\sH$, P.E.T.~Jorgensen \cite{J} and S.~\^Ota
\cite{Ot0}, \cite{Ot} have given a decomposition of a densely
defined operator $T$ as an operator sum of a \textit{closable}
operator, whose closure is again (the graph of) an operator, and a
singular operator, whose closure is the Cartesian product of a
closed subspace  of $\sH$ and a closed subspace  of $\sK$. This
decomposition is similar to a decomposition of nonnegative bounded
linear operators due to T.~Ando \cite{An} (see also \cite{Ko},
\cite{Ni}) and a decomposition of semibounded sesquilinear forms due
to B.~Simon (see \cite{Kosh}, \cite{S2}, \cite{S3}). It was pointed
out in these publications that there is an analogy with the Lebesgue
decomposition of a measure into a regular part and a singular part.

The purpose of this note is to show that there is a similar
decomposition in the case of linear relations from a Hilbert space
$\sH$ to a Hilbert space $\sK$. The notion of a linear relation as
a multivalued linear operator was introduced by R.~Arens \cite{A}
and extensively studied by E.A.~Coddington \cite{Co3} and by many
others. The treatment of Jorgensen and \^Ota for operators can be
relaxed: it is not necessary to consider operators which are
densely defined and in fact their treatment remains true for
relations. Indeed the language of relations seems to be the proper
context for such decompositions. Now the result is that any linear
relation from a Hilbert space $\sH$ to a Hilbert space $\sK$ has a
decomposition as an operator-like sum of a closable operator whose
closure is again (the graph of) an operator and a singular
relation whose closure is the Cartesian product of a closed
subspace of $\sH$ and a closed subspace of $\sK$. The components
of this decomposition can be characterized in various ways.

When the relation $T$ from $\sH$ to $\sK$ itself is considered
with the graph inner product, then its completion can be
contractively embedded into the Hilbert space $\sH$; the kernel of
this contraction corresponds to the multivalued part $\mul
\overline{T}$. This observation leads to a metric characterization
of the decomposition of the relation $T$. There is a similar
description for the decomposition of a pair of nonnegative
sesquilinear forms, cf. \cite{HSS}. By means of the above
mentioned result a metric characterization of closable operators
is presented.

The decomposition of Jorgensen and \^Ota also has connections with
the characteristic (projection) matrix introduced by
J.~von~Neumann \cite{vN} and M.H.~Stone \cite{Sto}, and Stone's
decomposition of a linear operator \cite{Sto}. The work of
A.E.~Nussbaum \cite{N} concerning orthogonal projections onto
closed subspaces of a Cartesian product is easily translated for
closed linear relations and this leads to the Stone decomposition
for closed linear relations, cf. \cite{M}. This makes it possible
to characterize the regular and singular parts of a linear
relation $T$ from a Hilbert space $\sH$ to a Hilbert space $\sK$
in terms of the orthogonal projection from the Cartesian product
$\sH \times \sK$ onto the closure $\overline{T}$.

It is also shown how the main decomposition result in this paper
(see Theorem~\ref{jor}) can be obtained by applying the general
characterization result of $\dom T^*$ in \cite{H} (cf. also
Lemma~\ref{lemma1} below); here $T^*$ is the adjoint relation of the
linear relation $T$.

\section{Preliminaries}

Here is a short review of notions associated with linear
relations. Recall that a linear relation $T$ from a Hilbert space
$\sH$ to a Hilbert space $\sK$ is a (not necessarily closed)
linear subspace of the Cartesian product $\sH \times \sK$. The
domain, range, kernel, and multivalued part of a linear relation
$T$ are defined by:
\[
  \dom T\okr\{\,f \in \sH:\, \{f,f'\} \in T\,\}, \quad \ran
T\okr\{\,f'\in \sK:\, \{f,f'\} \in T\,\},
\]
\[
  \ker T\okr\{\,f \in \sH :\, \{f,0\} \in T\,\}, \quad \mul
T\okr\{\,f' \in \sK:\, \{0,f'\} \in T\,\}.
\]
The formal inverse $T^{-1}$ is a linear relation from $\sK$ to
$\sH$ which is obtained from $T$ by interchanging the components
of the elements of $T$. Clearly, $\ran T=\dom T^{-1}$ and $\mul
T=\ker T^{-1}$. The linear relation $T$ is said to be closed if
$T$ is a closed subspace of the Cartesian product $\sH \times
\sK$. It $T$ is closed the kernel $\ker T$ and multivalued part
$\mul T$ of $T$ are automatically closed. The closures of the
domain and range of a linear relation $T$ are denoted by $\cdom T$
and $\cran T$. Observe that
\begin{equation}
\label{sing01} \dom \overline{T} \subset \cdom T\quad \mbox{and }
\quad \ran \overline{T} \subset \cran T.
\end{equation}
To see this, let $g \in \ran \overline{T}$, so that $\{f,g\} \in
\overline{T}$ for some $f \in \sH$. Then there exist elements
$\{f_n,g_n\} \in {T}$ such that $\{f_n, g_n\} \to \{f,g\}$. This
shows that $g \in \cran T$ and the second inclusion \eqref{sing01}
follows. The first identity is obtained from the second one by
inverting the relation $T$. The adjoint $T^*$ is a closed linear
relation from $\sK$ to $\sH$ defined by $T^*=JT^\perp=(JT)^\perp$,
where $J\{f,f'\}=\{f',-f\}$. Observe that $T^{**}=(T^\perp)^\perp$
is the closure $\overline{T}$ of $T$ in $\sH \times \sK$. Each of
the following identities is clear:
\begin{equation}
\label{einz+} (\ran T)^\perp=\ker T^*, \quad (\dom T)^\perp =\mul
T^*.
\end{equation}
Furthermore, the one identity in \eqref{einz+} is obtained by
inverting the relation in the other identity in \eqref{einz+}. When
$T$ in the identities in \eqref{einz+} is replaced by its adjoint
$T^*$, one also obtains
\begin{equation}
\label{zwei+} (\ran T^*)^\perp=\ker \overline{T}, \quad (\dom
T^*)^\perp =\mul \overline{T}.
\end{equation}
For two linear relations $A$ and $B$ from $\sH$ to $\sK$ there is
a \textit{componentwise sum} $A \hplus B$ from $\sH$ to $\sK$,
defined by
\[
  A \hplus B \okr \{\,\{f+g,f'+g'\} :\, \{f,f'\} \in A,\,\{g,g'\}
\in B\,\}.
\]
Note that\,\footnote{\;A dash ${^-}$ put aside still refers to the
closure operation.} $(A \hplus B)^*=A^* \cap B^*$ and $(A \hplus
B)^{-}= (\overline{A} \hplus \overline{B})^-$. The notation $A
\hoplus B$ is used to indicate that $A$ and $B$ are orthogonal in
the Cartesian product. As an example, observe that if $A$ is a
linear relation from $\sH$ to $\sK$ and $\sR$ is a linear subspace
of $\sK$, then the relation $T$ from $\sH$ to $\sK$, defined by
\begin{equation}
\label{summ}
  T=\{\,\{f,f'+\varphi\} :\, \{f,f'\} \in A, \, \varphi \in
\sR\,\},
\end{equation}
can be written as a componentwise sum $T=A \hplus B$, where the
linear relation $B$ from $\sH$ to $\sK$ is defined by $ B=
\{\,\{0,\varphi\} :\, \varphi \in \sR\, \}$. This leads to
\begin{equation}
\label{summ1}
 \overline{T}=(\overline{A} \hplus \{\,\{0,\varphi\} :\, \varphi
\in {\sR}^- \, \})^-.
\end{equation}
There is also an \textit{operator-like sum} $A +B$ from $\sH$ to
$\sK$, defined by
\[
  A + B \okr \{\,\{f,f'+f''\} :\, \{f,f'\} \in A,\,\{f,f''\} \in
B\,\},
\]
so that $\dom (A+B)=\dom A \cap \dom B$. The operator-like sum
reduces to the usual operator sum if $A$ and $B$ are (graphs of)
operators.

Now let $A$ be a linear relation from a Hilbert space $\sH$ to a
Hilbert space $\sR$, and let $B$ be a linear relation from the
Hilbert space $\sR$ to a Hilbert space $\sK$. The product $BA$ is
a linear relation from $\sH$ to $\sK$, defined by
\[
  BA \okr \{\,\{f,f'\} :\, \{f,g\} \in A, \, \{g,f'\} \in B\,\}.
\]
In general, $A^*B^* \subset (BA)^*$. However, when $B$ or $A^{-1}$
is (the graph of) a bounded everywhere defined operator on $\sR$,
then
\begin{equation}
\label{prod}
 (BA)^*=A^*B^*.
\end{equation}
This is known for the case of operators; for a proof in the case
of linear relations, see for instance \cite{HSS0}.

\section{Regular and singular relations}

A classical result of J.~von~Neumann states that the adjoint of a
densely defined linear operator is densely defined if and only if
the operator itself is closable (i.e., its closure is an
operator). The result which follows, and which can be deduced
immediately from the second identity in \eqref{zwei+}, does not
require the object in question to be {\sl a priori} (the graph of)
an operator; the latter comes out as an additional conclusion.
Thus call a linear relation $T$ from a Hilbert space $\sH$ to a
Hilbert space $\sK$ \textit{regular} if its closure $\overline{T}$
is (the graph of) an operator.

\begin{proposition}
\label{pro} Let $T$ be a linear relation from a Hilbert space
$\sH$ to a Hilbert space $\sK$. Then the following statements are
equivalent:
\begin{enumerate}
\def\labelenumi{\rm (\roman{enumi})}

\item $T$ is regular;

\item $\cdom T^*=\sK$,

\end{enumerate}
in which case $T$ is {\rm (}the graph of {\rm )} an operator.
\end{proposition}

Note that a linear operator $T$ is regular if and only if it is
closable. If a linear relation $T$ is regular, then $T$ is
automatically a closable operator. For any linear relation from a
Hilbert space $\sH$ to a Hilbert space $\sK$ the adjoint relation
$T^*$ from $\sK$ to $\sH$ is a closed linear relation; moreover,
the adjoint $T^*$ is a closed operator if and only if $\dom T$ is
dense in $\sH$ (independent of $T$ being regular). \\

A linear relation $T$ from a Hilbert space $\sH$ to a Hilbert
space $\sK$ is said to be \textit{singular} if
\begin{equation}
\label{sing}
  \ran T \subset \mul \overline{T} \quad \mbox{or equivalently}
\quad \cran T \subset \mul \overline{T}.
\end{equation}
The equivalence here is due to the closedness of $\mul
\overline{T}$. Furthermore, the inclusion
\begin{equation}
\label{sing0}
  \mul \overline{T} \subset \cran T,
\end{equation}
follows from \eqref{sing01} as $\mul \overline{T} \subset \ran
\overline{T}$. Therefore, a linear relation $T$ is singular if and
only if
\begin{equation}
\label{sing00}
  \cran T = \mul \overline{T},
\end{equation}
which follows by combining \eqref{sing} and \eqref{sing0}.

\begin{proposition}
\label{sin} Let $T$ be a linear relation from a Hilbert space
$\sH$ to a Hilbert space $\sK$. Then the following statements are
equivalent:
\begin{enumerate}
\def\labelenumi{\rm (\roman{enumi})}

\item $T$ is singular;

\item $\dom T^* \subset \ker T^*$ or, equivalently, $\dom T^*=\ker
T^*$;

\item $T^*=\dom T^* \times \mul T^*$;

\item $\overline{T}=\cdom T \times \mul \overline {T}$.

\end{enumerate}
\end{proposition}

\begin{proof}
(i) $\Rightarrow$ (ii) The identity in \eqref{sing00} implies
 that $(\cran T)^\perp=(\mul \overline{T})^\perp$, which is
equivalent to $\ker T^*=\cdom T^*$ by \eqref{einz+} and
\eqref{zwei+} . This implies that $\dom T^* \subset \ker T^*$,
while the reverse inclusion is obvious.

(ii) $\Rightarrow$ (iii) Let $\{f,g\} \in T^*$. Then in particular
$f \in \dom T^*$ and so $f \in \ker T^*$ by (ii). Therefore
$\{f,0\} \in T^*$ and this implies that $\{0,g\} \in T^*$, or $g
\in \mul T^*$. This shows that $\{f,g\} \in \dom T^* \times \mul
T^*$. Conversely, let $\{f,g\} \in \dom T^* \times \mul T^*$. Then
$\{0,g\} \in T^*$ by definition. Moreover, $f \in \dom T^*$ or by
(ii), $f \in \ker T^*$, so that $\{f,0\} \in T^*$. Clearly, by
linearity, this implies that $\{f,g\} \in T^*$.

(iii) $\Rightarrow$ (iv) Taking adjoints in (iii) yields
$T^{**}=(\mul T^*)^\perp \times (\dom T^*)^\perp$, which gives
(iv) by means of \eqref{einz+} and \eqref{zwei+}.

(iv) $\Rightarrow$ (i) It follows from $\ran \overline{T} = \mul
\overline{T}$ that $\ran T \subset \mul \overline{T}$. Hence, by
definition, $T$ is singular.
\end{proof}

\begin{corollary}
\label{sincor} Let $T$ be a linear relation from a Hilbert space
$\sH$ to a Hilbert space $\sK$. Then the following statements are
equivalent:
\begin{enumerate}
\def\labelenumi{\rm (\roman{enumi})}

\item $T$ is singular;

\item $T^{-1}$ is singular;

\item $T^*$ is singular.
\end{enumerate}
\end{corollary}

\begin{proof}
(i) $\Leftrightarrow$ (ii) Assume that $T$ is singular. Then part
(iv) of Proposition \ref{sin} implies that
\begin{equation}
\label{sid0} \ker\overline{T}=\cdom T \quad \mbox{or equivalently }
\mul\overline {T^{-1}}=\cran T^{-1}.
\end{equation}
Hence, $T^{-1}$ is singular by \eqref{sing00}. For the reverse
implication it is now enough to observe that $(T^{-1})^{-1}=T$.

(ii) $\Leftrightarrow$ (iii) Part (ii) of Proposition~\ref{sin}
shows that $T^*$ is singular if and only if $\dom \overline{T}=\ker
\overline{T}$. Clearly this is equivalent to \eqref{sid0}, which
means that $T^{-1}$ is singular.
\end{proof}

Note that if $T$ is singular, then in particular $\dom T^*=\ker
T^*$ is closed. A linear relation $T$ from a Hilbert space $\sH$
to a Hilbert space $\sK$ is said to be \textit{maximally singular}
if
\[
 \dom T^*=\{0\} \quad \mbox{or equivalently} \quad \mul
\overline{T}=\sK.
\]

\begin{proposition}
\label{ss} Let $T$ be a linear relation from a Hilbert space $\sH$
to a Hilbert space $\sK$. Then the following statements are
equivalent:
\begin{enumerate}
\def\labelenumi{\rm (\roman{enumi})}

\item $T$ is maximally singular;

\item $T$ is singular and $\cran T=\sK$;

\item $\overline{T}=\cdom T \times \sK$.

\end{enumerate}
\end{proposition}

\begin{proof}
(i) $\Rightarrow$ (ii) Assume that $T$ is maximally singular.
Clearly
\[
 \ran T \subset \sK=\mul \overline{T} \subset \ran \overline{T}
\subset \cran{T},
\]
where the last inclusion follows from \eqref{sing01}. Hence
\eqref{sing} is satisfied and $\cran T=\sK$.

(ii) $\Rightarrow$ (iii) Apply part (iv) of Proposition \ref{sin}
with $\mul \overline{T}=\cran T=\sK$.

(iii) $\Rightarrow$ (i) Clearly, $\mul \overline{T}=\sK$, so that
$T$ is maximally singular.
\end{proof}

As a simple example, observe that the linear relation $T$ defined
by $T=K^*U$ is maximally singular, when $U$ is any linear relation
and $K$ is an injective bounded linear operator with $\ran K
\subset (\dom U^*)^\perp$. To see this, note that $T^*=U^*K$, cf.
\eqref{prod}. Hence, if $\{f,f'\} \in T^*$, then $\{f,Kf\} \in K$
and $\{Kf, f'\} \in U^*$. This shows that $Kf \in \dom U^*$ which
leads to $Kf=0$. Since $K$ is injective, it follows that $f=0$.
Therefore $\dom T^*=\{0\}$ and $T$ is maximally singular, cf.
\cite{Ot}.

Observe that V.D.~Koshmanenko and S.~\^Ota \cite{KO} consider
(densely defined) linear operators $T$ which satisfy the property
$\dom T \subset \ker \overline{T}$. The inverse $T^{-1}$ of such an
operator is singular in the present sense, see \eqref{sing}, and
according to Corollary~\ref{sincor} then equivalently $T$ is
singular.

\section{Canonical decompositions}
\label{Sec4}

Let $T$ be a linear relation from a Hilbert space $\sH$ to a
Hilbert space $\sK$ and let $\overline{T}$ be its closure in the
Cartesian product $\sH \times \sK$. Denote the orthogonal
projection from $\sK$ onto $\mul \overline{T}$ by $P$. With $T$
are associated the linear relation $T_{\rm{reg}}$ from $\sH$ to
$\sK$ defined by
\begin{equation}
\label{oper1}
 T_{\rm{reg}}\okr\{\,\{f,(I-P)f'\} :\, \{f,f'\} \in T\,\},
\end{equation}
and the linear relation $T_{\rm{sing}}$ from $\sH$ to $\sK$
defined by
\begin{equation}
\label{oper2}
 T_{\rm{sing}}\okr\{\,\{f,Pf'\} :\, \{f,f'\} \in T\,\}.
\end{equation}
Observe that $T_{\rm{reg}}$ and $T_{\rm{sing}}$ have the same
domain $\dom T$. Moreover $T_{\rm{reg}}$ and $T_{\rm{sing}}$ are
(graphs of) operators if $T$ itself is (the graph of) an operator.
The following decomposition result is an adaptation of a result of
Jorgensen \cite{J}.

\begin{theorem} \label{jor}
Let $T$ be a linear relation from a Hilbert space $\sH$ to a
Hilbert space $\sK$. Then $T$ admits the canonical operator-like
sum decomposition
\begin{equation}
\label{can} T=T_{\rm{reg}} +T_{\rm{sing}} ,
\end{equation}
where $T_{\rm{reg}} $ is a regular relation from $\sH$ to $\sK$
and $T_{\rm{sing}} $ is a singular relation from $\sH$ to $\sK$
with
\[
\mul T_{\rm{sing}} =\mul T, \quad \mul ({T}_{\rm{sing}})^- =\mul
\overline{T}.
\]
\end{theorem}

\begin{proof}
To show that the relation $T_{\rm{reg}}$ in \eqref{oper1} is
regular, it suffices to show that its closure is an operator. Assume
therefore that there is a sequence $\{f_n,f_n'\} \in T$ such that
\[
 \{f_n, (I-P)f_n'\} \mapsto \{0,g\},
\]
which implies that $g \in (\mul \overline{T})^\perp$. Furthermore,
it follows from the definition of $T^*$ that for all $\{h,h'\} \in
T^*$
\begin{equation}
\label{nons}
\begin{split}
  0&=(f_n',h)-(f_n,h') \\
   &=((I-P)f_n',h)-(f_n,h')+(Pf_n',h) \\
   &=((I-P)f_n',h)-(f_n,h').
\end{split}
\end{equation}
Here the second identity in \eqref{zwei+} has been used. Taking
the limit $n \to \infty$ in \eqref{nons} leads to $(g,h)=0$ for
all $h \in \dom T^*$. Hence $g \in (\dom T^*)^\perp=\mul
\overline{T}$. Therefore $g=0$. It follows that $T_{\rm{reg}}$ is
a closable operator.

Next it will be shown that the relation $T_{\rm{sing}}$ is
singular. Note that $\{h,h'\} \in (T_{\rm{sing}})^*$ if and only
if
\[
 (h',f)=(h,Pf')=(P h,f') \quad \mbox{for all} \quad \{f,f'\} \in
T,
\]
or, equivalently, if and only if $\{P h, h'\} \in T^*$. Therefore
$h \in \dom (T_{\rm{sing}})^* $ if and only if $P h \in \dom T^*$.
Now observe that $P h \in \dom T^*$ if and only if $P h=0$, since
$\dom T^* \subset (\mul \overline{T})^\perp$. Furthermore, $Ph=0$
is equivalent to $h \in (\mul \overline{T})^\perp=\cdom T^*$.
Hence it follows that
\begin{equation}
\label{nons1}
  \dom (T_{\rm{sing}})^* =\cdom T^*.
\end{equation}
The same argument shows that $h \in \ker (T_{\rm{sing}})^*$ if and
only if $Ph \in \ker T^*$. Now, if $Ph \in \ker T^*$ then $Ph \in
\dom T^*$ and $Ph=0$. Conversely, if $Ph=0$ then $Ph \in \ker
T^*$. Hence, it follows that
\begin{equation}
\label{nons2}
  \ker (T_{\rm{sing}})^* =\cdom T^*.
\end{equation}
Combining \eqref{nons1} and \eqref{nons2} gives
\[
 \dom (T_{\rm{sing}})^* = \ker (T_{\rm{sing}})^*,
\]
in other words, the relation $T_{\rm{sing}}$ is singular by (ii)
of Proposition \eqref{sin}.

Moreover, \eqref{nons1} shows that
\[
  \mul (T_{\rm{sing}})^- = (\dom (T_{\rm{sing}})^*)^\perp =(\cdom
T^*)^\perp=\mul \overline{T}.
\]
Finally, observe that
\[
  \mul T_{\rm{sing}} =\{\,P f' :\, \{0,f'\}\in T\,\}= \{\,f' :\,
\{0,f'\}\in T\,\}=\mul T.
\]
This completes the proof.
\end{proof}

\begin{corollary}
The singular part of $T_{\rm{reg}}$ is the zero operator and the
regular part of $T_{\rm{sing}}$ is the zero operator. Hence
$T_{\rm{reg}}$ is equal to its regular part and $T_{\rm{sing}}$ is
equal to its singular part. In particular, $T$ is regular if and
only if $T=T_{\rm reg}$ and $T$ is singular if and only if
$T=T_{\rm sing}$.
\end{corollary}

\begin{proof}
The operator $T_{\rm{reg}}$ is closable; hence the multivalued part
of its closure is trivial. Therefore, its singular part is the zero
operator on $\dom T$.

The multivalued part of the closure of $T_{\rm{sing}}$ is equal to
$\mul (T_{\rm{sing}})^- =\mul \overline {T}$. Hence, the regular
part of $T_{\rm{sing}}$ is given by
\[
  \{\,\{f,(I-P)Pf'\} :\, \{f,f'\} \in T\,\},
\]
which is the zero operator on $\dom T$.

If $T=T_{\rm reg}$ ($T=T_{\rm sing}$), then $T$ is regular
(singular) by Theorem~\ref{jor}. Conversely, if $T$ is regular then
according to Proposition \ref{pro} $\cdom T^*=\sK$ or, equivalently,
$\mul \overline{T}=\{0\}$. Thus, $P=0$ and $T=T_{\rm reg}$.

Finally, if $T$ is singular then according to \eqref{sing00}
$Pf'=f'$ for all $f' \in \cran T$. In particular, $Pf'=f'$ for all
$f' \in \ran T$ and therefore $T=T_{\rm sing}$.
\end{proof}

\begin{corollary}
The singular part $T_{\rm{sing}}$ is maximally singular if and only
if $T$ is maximally singular.
\end{corollary}

\begin{proof}
Observe that the identity $\mul({T}_{\rm{sing}})^- =\mul
\overline{T}$ implies that $T_{\rm{sing}}$ is maximally singular if
and only if $T$ is maximally singular.
\end{proof}

The canonical decomposition in Theorem \ref{jor} is about the
decomposition of a linear relation $T$ as an operator-like sum of
a regular relation $T_{\rm{reg}}$ (a closable operator) and a
singular relation $T_{\rm{sing}}$. However, observe that because
\[
  \{f,f'\}=\{f, (I-P)f'+Pf'\}, \quad \{f,f'\} \in T,
\]
it can also be written as a component-wise sum
\[
  \{f,f'\}= \{f,(I-P)f'\} \hplus \{0,Pf'\}, \quad \{f,f'\}\in T
\]
which in fact is an orthogonal sum $\hoplus$ in $\sH\times\sK$.
This implies that
\begin{equation}
\label{deco}
  T \subset T_{{\rm reg}} \hoplus \{\,\{0,g\} :\, g \in \mul
\overline{T}\,\}.
\end{equation}
When $T$ is closed, this argument leads to the usual decomposition
which goes back to \cite{A} and \cite{Co3}.

\begin{proposition}
\label{ca} Let $T$ be a closed linear relation from a Hilbert
space $\sH$ to a Hilbert space $\sK$. Then $\mul T$ is closed,
$T_{\rm{reg}}$ is a closed operator, and
\begin{equation}
\label{deco1}
  T = T_{{\rm reg}} \hoplus \{\,\{0,g\} :\, g \in \mul {T}\,\}.
\end{equation}
\end{proposition}

\begin{proof}
If $T$ is closed, then $\mul \overline{T}=\mul T$ is closed and
$P$ is an orthogonal projection from $\sK$ onto $\mul T$. Hence in
\[
  \{f,f'\}=\{f, (I-P)f'+Pf'\}, \quad \{f,f'\} \in T ,
\]
one also has $\{0,Pf'\} \in T$, which leads to
\[
  \{f, (I-P)f'\} \in T \quad \mbox{for all} \quad \{f,f'\} \in T.
\]
It follows that $T_{{\rm reg}} \subset T$, and it is
straightforward to see that
\[
T_{{\rm reg}}=\{\,\{f,g\} \in T :\, g \perp \mul T\,\},
\]
which is clearly closed. Indeed, the righthand side of
\eqref{deco1} is contained in the lefthand side. The reverse
inclusion follows from \eqref{deco} keeping in mind that $\mul
\overline{T}=\mul T$.
\end{proof}

For similar orthogonal operator parts under the weaker condition
that only $\mul T$ is closed, see \cite{HSSz}.

\begin{proposition}
\label{over}
 Let $T$ be a linear relation from a Hilbert space $\sH$ to a
Hilbert space $\sK$. Then $(\overline{T})_{\rm reg}$ is closed and
\begin{equation}
\label{A}
 (\overline{T})_{\rm reg}=(T_{\rm reg})^-.
\end{equation}
Furthermore,
\begin{equation}
\label{B}
 ((\overline{T})_{\rm sing})^-=(T_{\rm sing})^-.
\end{equation}
\end{proposition}

\begin{proof}
The definition of regular part implies that ${T}_{\rm reg} \subset
(\overline{T})_{\rm reg}$. To see this let $\{f, (I-P)f'\}$ with
$\{f,f'\} \in T$ be an element in ${T}_{\rm reg}$. But then also
$\{f,f'\} \in \overline{T}$ and since $P$ is an orthogonal
projection onto $\mul \overline{T}$, it follows that $\{f,
(I-P)f'\}$ belongs to $(\overline{T})_{\rm reg}$.

The definition of regular part also implies that
$(\overline{T})_{\rm reg} \subset (T_{\rm reg})^-$. To see this
let $\{f, (I-P)f'\}$ with $\{f,f'\} \in \overline{T}$ be an
element of $(\overline{T})_{\rm reg}$. Then there exists a
sequence $\{f_n,f_n'\} \in T$ such that $\{f_n,f_n'\} \to
\{f,f'\}$. However this implies that the sequence
$\{f_n,(I-P)f_n'\} \in T_{\rm reg}$ approximates $\{f, (I-P)f'\}$.
In other words $(\overline{T})_{\rm reg} \subset (T_{\rm reg})^-$.

Combining these two assertions it follows that
\[
  {T}_{\rm reg} \subset (\overline{T})_{\rm reg} \subset (T_{\rm
reg})^-,
\]
which leads to
\[
  ((\overline{T})_{\rm reg})^-=(T_{\rm reg})^-.
\]
A similar argument for the singular part gives
\[
  ((\overline{T})_{\rm sing})^-=(T_{\rm sing})^-.
\]
Therefore \eqref{B} follows. Moreover, since $\overline{T}$ is
closed, it follows from Proposition \ref{ca} that
$\overline{T}_{\rm reg}$ is closed, which leads to \eqref{A}.
\end{proof}

\begin{corollary}
\label{nieu} Let $T$ be a linear relation from a Hilbert space
$\sH$ to a Hilbert space $\sK$. Then
\begin{equation}
\label{nieuw}
   (\overline{T}_{{\rm sing}})^-= \cdom T \times \mul
\overline{T}.
\end{equation}
\end{corollary}

\begin{proof}
By Proposition \ref{over} the lefthand side is equal to $(T_{\rm
sing})^-$. Moreover,
\[
 (T_{\rm sing})^-=\cdom T_{\rm sing} \times \mul (T_{\rm sing})^-
 =\cdom T \times \mul \overline{T},
\]
due to Proposition \ref{sin} and Theorem \ref{jor}.
\end{proof}

\begin{remark}
The relation $\overline{T}_{{\rm reg}}$ is closed; however, in
general, the relation $\overline{T}_{{\rm sing}}$ is not closed.
For assume that $\overline{T}_{{\rm sing}}$ is closed. Then by
Corollary \ref{nieu}
\[
  \overline{T}_{{\rm sing}}=\cdom \overline{T} \times \mul
\overline{T}.
\]
This implies that $\dom \overline{T}_{{\rm sing}}=\cdom
\overline{T}$, whereas by definition $\dom \overline{T}_{{\rm
sing}}=\dom \overline{T}$. Hence, if $\overline{T}_{{\rm sing}}$
is closed, then $\dom \overline{T}$ is closed.
\end{remark}

\section{A metric characterization}

Let $T$ be a linear relation from a Hilbert space $\sH$ to a
Hilbert space $\sK$, and introduce the graph inner product on $T$
by:
\[
 (\{f,f'\}, \{g,g'\})_T \okr (f,g)+(f',g'), \quad \{f,f'\},
 \{g,g'\} \in T.
\]
Denote the corresponding inner product space by $\sG(T)$. Define
the mapping $\iota_T$ by
\[
 \iota_T \{f,f'\} \okr f, \quad \{f,f'\} \in T,
\]
so that $\iota_T$ is a contraction on $\sG(T) $ with values in
$\dom T \subset \sH$. The isometric part of $\iota_T$ is its
restriction to $\ker T \times \{0\}$. Note that the Hilbert space
completion of $\sG(T)$ is $\sG(\overline{T})$, where
$\sG(\overline{T})$ stands for the construction related to the
closure $\overline{T}$ of $T$. Clearly $\sG(\overline{T}) =\sG(T)
$ if and only if the relation $T$ is closed, in which case the
contraction $\iota_T$ is closed. In general, the contraction
$\iota_T$ has a closure $\overline{\iota}_T$, which is a
contraction defined on all of $\overline{T}$. The contraction
$\iota_{\overline{T}}$ is also defined everywhere on
$\overline{T}$. Clearly the restrictions of the contractions
$\overline{\iota}_T$ and $\iota_{\overline{T}}$ coincide on $T$,
which is a dense set in $\overline{T}$; hence they coincide:
\[
   \overline{\iota}_T = \iota_{\overline{T}}.
\]
It follows from $\dom T \subset \dom \overline{T} \subset \cdom T$
that $\overline{\iota}_T = \iota_{\overline{T}}$ is a contraction
defined on all of $\sG(\overline{T}) $ with values in $\cdom T$.
Since $\iota_{\overline{T}}$ maps $\sG(\overline{T})$ into $\cdom
T$, it follows by the first identity in \eqref{einz+} that
\[
  \cran (\iota_T)^*=\cran (\iota_{\overline{T}})^*=
\sG(\overline{T}) \ominus \ker \iota_{\overline{T}}.
\]
In particular, this shows that
\begin{equation}
\label{range}
  (\iota_T)^* \uphar \dom T \quad \mbox{is dense in} \quad
\sG(\overline{T}) \ominus \ker \iota_{\overline{T}}.
\end{equation}
The above construction is the analog for relations of a
construction involving semibounded sesquilinear forms, cf.
\cite{HSS}, \cite{S3}. The following lemma is obvious.

\begin{lemma}
\label{help} Let $T$ be a linear relation from a Hilbert space
$\sH$ to a Hilbert space $\sK$. Then
\begin{equation}
\label{help1}
  \ker \iota_{\overline{T}}= \{0\} \times
\mul \overline{T}.
\end{equation}
\end{lemma}

Denote by $Q $ the orthogonal projection from the Hilbert space
$\sG(\overline{T})$ onto its closed linear subspace $\ker
\iota_{\overline{T}} $.

\begin{lemma}
\label{help3}
 Let $T$ be a linear relation from a Hilbert space
$\sH$ to a Hilbert space $\sK$. The orthogonal projection $Q $ in
$\sG(\overline{T})$ and the orthogonal projection $P$ in $\sK$ are
related by
\[
 Q \{\varphi,\varphi'\}=\{0,P\varphi'\}, \quad (I-Q
)\{\varphi,\varphi'\}=\{\varphi,(I-P)\varphi'\}, \quad
\{\varphi,\varphi'\} \in T.
\]
\end{lemma}

\begin{proof}
In the sense of the inner product of the Hilbert space
$\sG(\overline{T})$ each element $\{\varphi,\varphi'\} \in T$ has
the orthogonal decomposition
\begin{equation}
\label{uksi} \{\varphi,\varphi'\}= \{h,h'\} \hplus \{k,k'\},
\end{equation}
where
\begin{equation}
\label{uksi0}
\{h,h'\}=(I-Q )\{\varphi,\varphi'\}, \quad
\{k,k'\}=Q \{\varphi,\varphi'\}.
\end{equation}
The elements $\{h,h'\}$ and $\{k,k'\}$ are orthogonal in the sense
of $\sG(\overline{T})$:
\[
  0=(h,k)+(h',k').
\]
It follows from Lemma \ref{help} that $k=0$, so that $h=\varphi$
and $(h',k')=0$. Hence the decomposition \eqref{uksi} can also be
written as
\begin{equation}
\label{kaksi} \{\varphi,\varphi'\}= \{\varphi,h'+k'\}, \quad k'
\in \mul \overline{T},
\end{equation}
where the decomposition of $\varphi'=h'+k'$ is orthogonal in
$\sK$. In other words
\[
  h' =(I-P) \varphi', \quad k'=P \varphi',
\]
so that \eqref{uksi0} reads as
\[
  (I-Q )\{\varphi,\varphi'\}=\{\varphi,(I-P) \varphi'\}, \quad
  Q\{\varphi,\varphi'\}= \{0,P \varphi'\},
\]
which gives the statement of the lemma.
\end{proof}

The construction involving the contraction $\iota_T$ now leads to
a metric characterization of the elements in the ranges of $T_{\rm
reg}$ and $T_{\rm sing}$.

\begin{lemma}
\label{three}
Let $T$ be a linear relation from the Hilbert space
$\sH$ to the Hilbert space $\sK$. Then for all
$\{\varphi,\varphi'\} \in T$:
\begin{equation}
\label{drei} \|P\varphi'\|^2=
 \|\varphi'\|^2+ \inf_{h \in \dom T} \left\{\|\varphi +h\|^2-
 \inf_{\{g,g'\} \in T} \left\{\|g'\|^2+\|g+h\|^2 \right\}
 \right\},
\end{equation}
and
\begin{equation}
\label{vier} \|(I-P)\varphi'\|^2=\sup_{h\in \dom T} \inf_{\{g,g'\}
\in T} \left\{\, \|g+h\|^2-\|\varphi+h\|^2+\|g'\|^2 \,\right\}
\end{equation}
\end{lemma}

\begin{proof}
According to Lemma \ref{help3} $Q
\{\varphi,\varphi'\}=\{0,P\varphi'\}$, which implies that
\begin{equation}
\label{vijf}
 \|P\varphi'\|^2= \|Q \{\varphi,\varphi'\}\|_{\overline{T}}^2.
\end{equation}
Recall that $\sG(\overline{T})=\ker \iota_{\overline{T}} \oplus
(\sG(\overline{T}) \ominus \ker \iota_{\overline{T}})$, which
gives
\[
 \|Q \{\varphi,\varphi'\}\|_{\overline{T}}^2= \inf \left\{\,
\|\{\varphi,\varphi'\}-\{\alpha,\alpha'\}\|_{\overline{T}}^2:\,
\{\alpha,\alpha'\} \in \sG(\overline{T}) \ominus \ker
\iota_{\overline{T}} \,\right\}.
\]
However, since $ \iota_T^*(\dom T)$ is dense in $\sG(\overline{T})
\ominus \ker \iota_{\overline{T}} $, it follows that
\begin{equation}
\label{einz}
\begin{split}
 \left\|Q \{\varphi,\varphi'\}\right\|^2_{\overline{T}} &= \inf_{h
 \in \dom T} \left\{\left(\{\varphi,\varphi'\}+\iota_T^* h,
\{\varphi,\varphi'\}+\iota_T^* h \right)_{\overline{T}} \right\}   \\
&= \inf_{h \in \dom T}
\left\{(\varphi,\varphi)+(\varphi',\varphi')+(\varphi,h)+(h,\varphi)
+\left( \iota_T^* h, \iota_T^*h\right)_{\overline{T}} \right\}
\\
&=\|\varphi'\|^2+\inf_{h \in \dom T} \left\{ \|\varphi
+h\|^2-\|h\|^2 +\left\|\iota_T^* h \right\|^2_{\overline{T}}
\right\}.
\end{split}
\end{equation}
Furthermore, since $T$ is dense in $\sG(\overline{T})$, every
element of the form $\iota_T^* h$, $h \in \dom T$, can be
approximated by elements in $T$, which leads to
\begin{equation}
\label{zwei}
\begin{split}
 0&= \inf_{\{g,g'\} \in T}
 \left\{ \left\|\{g,g'\}+\iota_T^* h \right\|^2_{\overline{T}}  \right\} \\
  &= \left\| \iota_T^* h \right\|^2_{\overline{T}}+\inf_{\{g,g'\}
\in T} \left\{
(g,g)+(g',g')+(g,h)+(h,g) \right\}     \\
  &= -\|h\|^2 + \left\| \iota_T^* h \right\|^2_{\overline{T}}
+\inf_{\{g,g'\} \in T}
  \left\{\|g'\|^2+\|g+h\|^2\right\}.
\end{split}
\end{equation}
Combining the identities \eqref{einz}, \eqref{zwei}, and
\eqref{vijf} gives \eqref{drei}. Clearly \eqref{vier} follows from
\eqref{drei}.
\end{proof}

The above lemma leads to a metric characterization of closable
operators.

\begin{theorem}
\label{deux} Let $T$ be a linear operator from a Hilbert space
$\sH$ to a Hilbert space $\sK$. Then $T$ is closable if and only
if for all $\varphi \in \dom T$:
\begin{equation}
\label{een1}
  \|T\varphi\|^2=\sup_{h\in \dom T} \inf_{g \in \dom T} \left\{\,
\|g+h\|^2-\|\varphi+h\|^2+\|Tg\|^2 \,\right\}.
\end{equation}
\end{theorem}

\begin{proof}
Assume that \eqref{een1} holds for all $\varphi \in \dom T$. By
Lemma \ref{three} this means that
\[
  \|T\varphi\|^2=\|(I-P)T\varphi\|^2, \quad \varphi \in \dom T.
\]
In other words, $PT\varphi=0$ and hence $T\varphi \in \ker P$ for
all $\varphi \in \dom T$. Therefore $\ran T \subset \ker P$, which
leads to
\[
  \ran \overline{T} \subset \ker P, \quad \mbox{in particular},
\quad \mul \overline{T} \subset \ker P.
\]
Since $P$ is an orthogonal projection onto $\mul \overline{T}$,
this implies that $\mul \overline{T}=\{0\}$, i.e., $\overline{T}$
is an operator. Hence, the operator $T$ is closable.

Conversely, assume that the operator $T$ is closable. Then according
to Theorem \ref{jor} $T=T_{\rm{reg}}$ and $P=0$. Hence if
$\{\varphi,\varphi'\} \in T$, then $T\varphi=\varphi'=(I-P)\varphi'$
and the result follows from \eqref{vier}.
\end{proof}

In general the supremum and the infimum in \eqref{een1} are not
attained. However, when the operator $T$ is densely defined and
closed one can say more.

\begin{lemma}
\label{Tclosed} Let $T$ be a densely defined closed linear
operator from $ \sH$ to $\sK$. Then
\begin{equation}
\label{las1} \min_{g \in \dom T} \left(\|g+h\|^2 +\|Tg\|^2
\right)=\|h\|^2-\|(I+T^*T)^{-\half}h\|^2, \quad h \in \sH,
\end{equation}
and the minimum is attained for $g=-(I+T^*T)^{-1}h \in \dom T^*T$.
\end{lemma}

\begin{proof}
Since $T$ is densely defined and closed, $T^*T$ is a selfadjoint
operator, which is nonnegative. Observe that for all $g \in \dom
T$ there is the identity
\[
 \|g\|^2+\|Tg\|^2 =\|(I+T^*T)^\half g\|^2.
\]
Hence, for all $g\in \dom T$ and $h \in \sH$ one has
\[
\begin{split}
  \|g+h\|^2 &+\|Tg\|^2 \\&=\|h\|^2-\|(I+T^*T)^{-\half}h\|^2
  +\left\|(I+T^*T)^\half g+(I+T^*T)^{-\half} h\right\|^2
  \\
  & \ge \|h\|^2-\|(I+T^*T)^{-\half}h\|^2.
\end{split}
\]
Due to $\ran (I+T^*T)^\half=\sH$ this implies that \eqref{las1}
holds and the minimum is attained for $g=-(I+T^*T)^{-1}h \in \dom
T^*T $.
\end{proof}

A combination of Theorem \ref{deux} and Lemma \ref{Tclosed} leads
to the following characterization.

\begin{proposition}
Let $T$ be a densely defined closed linear operator from $ \sH$ to
$\sK$. Then
\begin{equation}
\label{een1++}
  \|T\varphi\|^2=\sup_{h\in \dom T} \min_{g \in \dom T} \left\{\,
\|g+h\|^2-\|\varphi+h\|^2+\|Tg\|^2 \,\right\},
\end{equation}
and the supremum is a maximum if and only if $\varphi \in \dom
T^*T$.
\end{proposition}

\begin{proof}
Since the operator $T$ is closed the identity \eqref{een1} holds
for all $\varphi \in \dom T$. Clearly, it follows from
\eqref{een1} and \eqref{las1} that
\begin{equation}
\label{een1+} \|T\varphi\|^2=\sup_{h\in \dom T} \left(\|h\|^2
 -\|(I+T^*T)^{-\half}h\|^2 -\|\varphi+h\|^2 \right), \quad \varphi
\in \dom T .
\end{equation}
Observe that for all $\varphi \in \dom T$ and $h \in \sH$ one
has
\[
\begin{split}
   \|h\|^2 &-\|(I+T^*T)^{-\half}h\|^2 -\|\varphi+h\|^2 \\
&= -\|(I+T^*T)^{-\half} h +(I+T^*T)^{
\half}\varphi\|^2+\|T\varphi\|^2 ,
\end{split}
\]
which implies that the supremum in \eqref{een1+} and hence in
\eqref{een1++} is a maximum if and only if $\varphi \in \dom
T^*T$.
\end{proof}

In particular, if $T$ is a bounded linear operator from $\sH$ to
$\sK$ then the supremum in \eqref{een1++} can be replaced by a
maximum. The original observation about the minimum and maximum in
\eqref{een1} for bounded linear operators goes back to
L.~L\'aszl\'o (personal communication, 2005).

\begin{corollary}
\label{cordeux} Let $T$ be a closed linear operator from a Hilbert
space $\sH$ to a Hilbert space $\sK$ and let $S$ be a linear
operator from the Hilbert space $\sH$ to the Hilbert space $\sR$
with $\dom S=\dom T$ and
\begin{equation} \label{trois}
  \|T\varphi\| \le \|S\varphi\|, \quad \varphi \in \sD\okr\dom S=\dom T.
\end{equation}
Then for all $\varphi \in \sD$:
\begin{equation}
\label{een2} \|T\varphi\|^2 \le \sup_{h\in \sD} \inf_{g \in \sD}
\left\{\, \|g+h\|^2-\|\varphi+h\|^2+\|Sg\|^2 \,\right\} \le
\|S\varphi\|^2.
\end{equation}
\end{corollary}

\begin{proof}
It follows from Theorem \ref{deux} and the assumption \eqref{trois}
that
\[
\begin{split}
  \|T\varphi\|^2&=\sup_{h\in \sD} \inf_{g \in \sD} \left\{\,
\|g+h\|^2-\|\varphi+h\|^2+\|Tg\|^2 \,\right\} \\
&\le \sup_{h\in \sD} \inf_{g \in \sD} \left\{\,
\|g+h\|^2-\|\varphi+h\|^2+\|Sg\|^2 \,\right\}, \quad \varphi \in
\sD.
\end{split}
\]
The choice $g=\varphi$ leads to the inequality
\[
 \sup_{h\in \sD} \inf_{g \in \sD} \left\{\,
 \|g+h\|^2-\|\varphi+h\|^2+\|Sg\|^2 \,\right\} \le \|S\varphi\|^2,
\]
which results in \eqref{een2}.
\end{proof}

The result in Theorem~\ref{deux} can be interpreted in terms of
parallel sums and differences (see \cite{HSS}), cf. \cite{An},
\cite{AnS}, \cite{Ni}, \cite{Szy}. Furthermore, by replacing $T$
with its inverse $T^{-1}$ one obtains a similar characterization for
the implication $\ker T=\{0\}\Rightarrow \ker \overline{T}=\{0\}$.
Observe also that if the operator $T$ in Theorem~\ref{deux} is not
closable, then for some $\varphi\in\dom T$:
\[
  \|T\varphi\|^2 > \sup_{h\in \dom T} \inf_{g \in \dom T} \left\{\,
\|g+h\|^2-\|\varphi+h\|^2+\|Tg\|^2 \,\right\}.
\]

\section{The Stone decomposition}

Since a closed linear relation from a Hilbert space $\sH$ to a
Hilbert space $\sK$ is by definition a closed linear subspace of
the Cartesian product $\sH \times \sK$, there is a one-to-one
correspondence between the orthogonal projections in $\sH \times
\sK$ and the closed linear relations from $\sH$ to $\sK$. This
section gives a short review of the consequences of this
correspondence, cf. \cite{N}, \cite{Sto}.

Let $\sH$ and $\sK$ be Hilbert spaces and let $R$ be an orthogonal
projection on the Cartesian product $\sH \times \sK$. Decompose
$R$ according to the Cartesian product
\begin{equation}
\label{proj}
  R=\begin{pmatrix} R_{11} & R_{12} \\ R_{21} & R_{22}
    \end{pmatrix},
\end{equation}
so that $R_{11} \ge 0$, $R_{22} \ge 0$, and $R_{21}=R_{12}^*$.

\begin{lemma} [\cite{N}]
\label{nuss} The entries in the block decomposition \eqref{proj}
satisfy:
\begin{enumerate}
\def\labelenumi{\rm (\roman{enumi})}
\item $\ker R_{11} \subset \ker R_{21}$;

\item $\ker R_{22} \subset \ker R_{12}$;

\item $\ker (I-R_{11}) \subset \ker R_{21}$;

\item $\ker (I-R_{22}) \subset \ker R_{12}$.

\end{enumerate}
\end{lemma}

Let $T$ be a closed relation from $\sH$ to $\sK$ and let $R$ be
the corresponding orthogonal projection onto $T$. The
corresponding matrix $R=(R_{ij})$ is called the
\textit{characteristic matrix} of $T$. It follows by definition
that
\begin{equation}
\label{t}
   {T}=\{\,\{ R_{11}h+ R_{12}h', R_{21}h+ R_{22}h' \} :\, \{h,h'\}
  \in \sH \times \sK\,\}.
\end{equation}
The representation \eqref{t} is called the \textit{Stone
decomposition} of $T$. It is clear that
\begin{equation} \label{proj1}
   \begin{pmatrix} R_{22} & R_{21} \\ R_{12} & R_{11}
    \end{pmatrix}
\end{equation}
is the characteristic matrix of the inverse relation $T^{-1}$
(keeping the order of the parametrizing elements in mind).

\begin{lemma}
\label{nus1} Let $T$ be a closed linear relation from a Hilbert
space $\sH$ to a Hilbert space $\sK$, and let $R $ be its
characteristic matrix. Then
\begin{equation}
\label{domran}
  \dom T=\ran \begin{pmatrix} R_{11} & R_{12} \end{pmatrix}, \quad
  \ran T=\ran \begin{pmatrix} R_{21} & R_{22} \end{pmatrix},
\end{equation}
and
\begin{equation}
\label{domran1}
  \ker T=\ker (I-R_{11}), \quad \mul T=\ker (I-R_{22}).
\end{equation}
\end{lemma}

\begin{proof}
The identities in \eqref{domran} are clear. Observe also that
$\ran T=\dom T^{-1}$, so that, in fact, the second identity in
\eqref{domran} follows from the first identity and \eqref{proj1}.

To prove the first identity in \eqref{domran1}, let $h \in \ker
T$, then $\{h,0\} \in T$, so that in particular
$\{h,0\}=R\{h,0\}$, or equivalently,
\[
 \begin{pmatrix} h \\ 0 \end{pmatrix}=
 \begin{pmatrix} R_{11} & R_{12} \\ R_{21} & R_{22}
    \end{pmatrix} \begin{pmatrix} h \\ 0 \end{pmatrix}
    =\begin{pmatrix} R_{11}h \\ R_{21}h \end{pmatrix},
\]
which implies that $R_{11}h=h$ and $R_{21}h=0$, and, by Lemma
\ref{nuss}, this is equivalent to $h \in \ker (I-R_{11})$.

The statement concerning $\mul T=\ker T^{-1}$ follows from the
first identity in \eqref{domran1} and \eqref{proj1}.
\end{proof}

\begin{corollary}
Let $T$ be a closed linear relation from a Hilbert space $\sH$ to
a Hilbert space $\sK$. Its regular part $T_{\rm{reg}}$ is given by
\begin{equation}
\label{t0} T_{\rm{reg}} =\{\,\{ R_{11}h+ R_{12}k, R_{21}h+ R_{22}k
\} :\, h \in \sH, \, k \in (\mul T)^\perp \,\}.
\end{equation}
\end{corollary}

\begin{proof}
Consider the representation \eqref{t} of the closed linear
relation $T$. Decompose the variable $h' \in \sH$ by $h'=k +
\varphi$ with $k \in (\mul T)^\perp$ and $\varphi \in \mul T$.
Then it follows from Lemma \ref{nuss} and Lemma \ref{nus1} that
\[
  R_{11}h+ R_{12}h'=R_{11}h+ R_{12}k, \quad R_{21}h+
  R_{22}h'=R_{21}h+ R_{22}k+ \varphi,
\]
and it is also clear that $R_{21}h+ R_{22}k \in (\mul T)^\perp$.
This completes the proof.
\end{proof}

The definition of the adjoint $T^{*}=JT^\perp$ (where the product
in the righthand side is carried out in the indicated order) gives
the corresponding characteristic matrix
\begin{equation}
\label{proj2}
   \begin{pmatrix} I-R_{22} & R_{21} \\ R_{12} & I-R_{11}
    \end{pmatrix},
\end{equation}
which leads to the following parametrization of $T^*$:
\begin{equation}
\label{t1}
   T^*=\{\,\{ (I-R_{22})h+ R_{21}h', R_{12}h+ (I-R_{11})h' \} :\,
  \{h,h'\} \in \sH \times \sK\,\}.
\end{equation}
By Lemma \ref{nus1} it follows that
\begin{equation}
\label{domran2+}
  \dom T^*= \ran \begin{pmatrix} (I-R_{22}) & R_{21}
\end{pmatrix},
\quad \ran T^*=\ran \begin{pmatrix} R_{12} & (I-R_{11})
\end{pmatrix},
\end{equation}
and
\begin{equation}
\label{domran2}
 \ker T^*= \ker R_{22}, \quad \mul T^* =\ker R_{11}.
\end{equation}
The regular part of $T^*$ is given by
\[
  (T^*)_{{\rm reg}}=\{\,\{ (I-R_{22})h+ R_{21}k, R_{12}h+
  (I-R_{11})k \} :\, h \in \sH, \, k\in (\mul T^*)^\perp \,\}.
\]
The following result is now straightforward, see \cite{Sto} for
the case of operators and \cite{M} for the case of relations.

\begin{lemma}
Let $T$ be a closed linear relation from a Hilbert space $\sH$ to
a Hilbert space $\sK$. Then the products $T^*T$ and $TT^*$ are
nonnegative selfadjoint relations in $\sH$ and $\sK$,
respectively, and
\[
  (T^*T+I)^{-1} \quad \mbox{and} \quad (TT^*+I)^{-1}
\]
are {\rm (}the graphs of {\rm )} bounded linear operators defined
on all of $\sH$ and $\sK$, respectively. The characteristic matrix
$R$ of $T$ is given by
\[
  R=\begin{pmatrix}(T^* {T}+I)^{-1}& (T^*)_{\rm{reg}}(TT^*+I)^{-1} \\
{T}_{\rm{reg}} (T^* {T}+I)^{-1} & I-(TT^*+I)^{-1}
\end{pmatrix}.
\]
\end{lemma}

\begin{proof}
Let $h \in \sH$, then there is a unique decomposition
\[
 \{h,0\} =\{\varphi,\varphi'\}+\{\psi,\psi'\}, \quad
 \{\varphi,\varphi'\} \in {T}, \quad \{\psi,\psi'\} \in
 T^\perp=JT^*,
\]
since $\sH \times \sK= {T} \hoplus T^\perp$. Hence
\[
 h=\varphi+\psi, \quad \varphi'+\psi'=0,
\]
which leads to $\{\psi,\psi'\}=\{\psi,-\varphi'\}\in JT^*$ and
$\{\varphi',\psi\} \in T^*$. Therefore, $\{\varphi,\psi\} \in T^*
{T}$ and
\[
 \{\varphi,h\}=\{\varphi,\varphi+\psi\} \in T^* {T}+I,
\]
so that $h \in \ran (T^* {T}+I)$. Thus $\ran (T^* {T}+I)=\sH$ and
$\varphi=(T^* {T}+I)^{-1}h$. Note that $\{\psi,\psi'\} \in
T^\perp=JT^*$ or, equivalently, $\{\psi',-\psi\} \in T^*$, so that
$\varphi'=-\psi'\in \dom T^* \subset (\mul {T})^\perp$ and it
follows that $\varphi'= {T}_{\rm{reg}}\varphi$. Therefore
\[
\varphi=(T^* {T}+I)^{-1}h, \quad \varphi'= {T}_{\rm{reg}} (T^*
{T}+I)^{-1}h.
\]
Hence the first column of $R $ is completely determined. By
formally replacing $T$ by $T^*$ also the first column of the
characteristic matrix of $T^*$ is determined. Now the second
column of $R$ is obtained via \eqref{proj2}.
\end{proof}

The statements concerning the products $T^*T$ and $TT^*$ are known
(cf. \cite{HSSW1}), but the proof is repeated here to make the
identification of the first column of $R$ understandable.

\section{A characterization of the regular and singular parts
via the Stone decomposition}

Let $T$ be a linear relation from a Hilbert space $\sH$ to a
Hilbert space $\sK$ and let $\overline{T}$ be its closure.
Decompose the Cartesian product $\sH \times \sK$ as follows
\begin{equation}
\label{stone}
      \sH \times \sK =\overline{T} \hoplus {T}^\perp,
\end{equation}
where $\hoplus$ denotes the orthogonal component-wise sum. Let
$\widetilde{R}$ denote the orthogonal projection from $\sH \times
\sK$ onto $\overline{T}$ and decompose $\widetilde{R}$ according
to the Cartesian product:
\[
\widetilde{R}=
\begin{pmatrix} \widetilde{R}_{11} & \widetilde{R}_{12} \\
\widetilde{R}_{21} & \widetilde{R}_{22}
\end{pmatrix}.
\]
Observe that
\begin{equation}
\label{stone0}
  \mul \overline{T}=\ker (I-\widetilde{R}_{22}) \subset \ker
\widetilde{R}_{12},
\end{equation}
by combining Lemma \ref{nus1} and Lemma \ref{nuss}, and that
\begin{equation}
\label{stone00}
  \mul T^*=\ker \widetilde{R}_{11} \subset \ker
\widetilde{R}_{21},
\end{equation}
by combining \eqref{domran2} and Lemma \ref{nuss}. In particular,
this implies that
\begin{equation}
\label{stone01}
   \widetilde{R}_{11} \uphar_{\mul T^*}=0\uphar_{\mul T^*}, \quad
\widetilde{R}_{21} \uphar_{\mul T^*}=0\uphar_{\mul T^*},
\end{equation}
and
\begin{equation}
\label{stone02}
   \widetilde{R}_{12} \uphar_{\mul \overline{T}}=0\uphar_{\mul
\overline{T}}, \quad \widetilde{R}_{22} \uphar_{\mul
\overline{T}}=I \uphar_{ \mul \overline{T}}.
\end{equation}

\begin{proposition}
The relation $T$ is regular if and only if
\begin{equation}
\label{stone3}
  \ker (I-\widetilde{R}_{22})=\{0\}.
\end{equation}
\end{proposition}

\begin{proof}
By definition the relation $T$ is regular if and only if
$\overline{T}$ is an operator, or, equivalently, $\mul
\overline{T}=\{0\}$. Hence the identity \eqref{stone3} follows
from the first identity in \eqref{stone0}.
\end{proof}

\begin{proposition}
\label{7.2} The relation $T$ is singular if and only if
\begin{equation}
\label{stone3+}
 \widetilde{R}_{12}=0\quad \mbox{or equivalently}
 \quad \widetilde{R}_{21}=0,
\end{equation}
in which case
\begin{equation}
\label{stone3++}
  \widetilde{R}=\begin{pmatrix} I \uphar_{ \cdom {T}} \oplus
0\uphar_{(\dom T)^\perp} & 0\\0 & I\uphar_{ \cran {T}} \oplus
0\uphar_{(\ran T)^\perp}\end{pmatrix}.
\end{equation}
\end{proposition}

\begin{proof}
Let $T$ be a singular relation. Then $\ran T \subset \mul
\overline{T}$ and it follows from \eqref{stone0} that for all
$\{\varphi,\varphi'\} \in T \subset \overline{T}$
\[
\begin{pmatrix} \varphi \\ \varphi' \end{pmatrix}
 =\begin{pmatrix} \widetilde{R}_{11} & \widetilde{R}_{12} \\
\widetilde{R}_{21} & \widetilde{R}_{22}
\end{pmatrix} \begin{pmatrix} \varphi \\ \varphi' \end{pmatrix}
= \begin{pmatrix} \widetilde{R}_{11} \varphi \\ \widetilde{R}_{21}
\varphi +\varphi'
\end{pmatrix},
\]
as $\varphi' \in \mul \overline{T}$. Hence for all $\varphi \in \dom
T$:
\begin{equation}
\label{newb2} \widetilde{R}_{11} \varphi=\varphi\quad \mbox{and}
\quad \widetilde{R}_{21} \varphi=0.
\end{equation}
Since the entries in the block decomposition of $\widetilde{R}$ are
bounded operators, one has $\widetilde{R}_{21}\uphar \cdom T=0\uphar
\cdom T$, which together with the second identity in \eqref{stone01}
shows that $\widetilde R_{21}=0$. This gives \eqref{stone3+}.

Conversely, assume that \eqref{stone3+} holds. Then necessarily
$\widetilde R_{22}$ is an orthogonal projection in $\sK$ and now
Lemma \ref{nus1} gives
\begin{equation}
\label{news1}
  \ran T \subset \ran \overline{T}=\ran \widetilde R_{22}=
  \ker (I-\widetilde R_{22})=\mul \overline{T},
\end{equation}
so that $T$ is singular.

It remains to discuss the matrix representation \eqref{stone3++}.
The representation of $\widetilde R_{11}$ is obtained from the first
identities in \eqref{stone01} and \eqref{newb2}. The representation
of $\widetilde R_{22}$ follows immediately from \eqref{news1}, since
$\widetilde R_{22}$ is an orthogonal projection in $\sK$.
\end{proof}

The result in Proposition~\ref{7.2} shows that there is a one-to-one
correspondence between the closed singular relations from $\sH$ to
$\sK$ and the orthogonal projections $R$ in $\sH\times\sK$ which are
of the form $R=R_{11}\oplus R_{22}$, where $R_{11}$ in an orthogonal
projection in $\sH$ and $R_{22}$ in an orthogonal projection in
$\sK$. Maximally singular relations can be described as follows.

\begin{corollary}
The relation $T$ is maximally singular if and only if
\begin{equation}
\label{stone4aa}
   \widetilde{R}_{22}=I,
\end{equation}
in which case \eqref{stone3+} holds and
\begin{equation}
\label{stone4}
   \widetilde{R}=\begin{pmatrix} I \uphar_{ \cdom {T}} \oplus
0\uphar_{(\dom T)^\perp} & 0\\0 & I \end{pmatrix}.
\end{equation}
\end{corollary}

\begin{proof}
By Lemma~\ref{nus1} (see \eqref{stone0})
$\mul\overline{T}=\ker(I-\widetilde R_{22})$ and thus $T$ is
maximally singular, i.e. $\mul \overline{T}=\sK$, if and only if
$\widetilde R_{22}=I$. Observe that \eqref{stone4aa} implies
\eqref{stone3+} in view of \eqref{stone0}.

The representation \eqref{stone4} is immediate from
\eqref{stone3++}.
\end{proof}

\section{Examples of canonical decompositions}

This section contains some illustrative examples concerning the
canonical decomposition of (not necessarily densely defined) linear
operators and relations.

\begin{example}
\label{StSz} \cite{StSz} Let $S$ be a closable operator in a
Hilbert space $\sH$ and let $Z=\span\{e,f\}$ be a one-dimensional
subspace in the Cartesian product $\sH \times \sH$. Note that
$\{e,f\} \in S$ if and only if $e \in \dom S$ and $Se=f$. Hence
the assumption $S \cap Z=\{0,0\}$ is equivalent to the assumption
that either $e \in \dom S$ and $Se \ne f$ or $e \notin \dom S$.
Under this assumption the component-wise sum
\begin{equation}
\label{defe} A=S \hplus Z,
\end{equation}
defines a one-dimensional extension (of the graph) of $S$. It
follows from \eqref{defe} that $A^{**}=S^{**} \hplus Z$, since $Z$
is one-dimensional. The one-dimensional extension $A$ is not regular
if and only if there is a nontrivial $k \in \sH$ for which $\{0,k\}
\in \overline{A}=A^{**}$, i.e.,
\[
 \{0,k\}=\{h, \overline{S}h\} \hplus \lambda \{e,f\},
\]
for some $h \in \dom \overline{S}=\dom S^{**}$ and $\lambda \in
\dC$. Equivalently, the one-dimensional extension $A$ is not
regular if and only if
\begin{equation}
\label{defe1} e \in \dom \overline {S}, \quad \overline {S}e \ne
f,
\end{equation}
in which case
\[
  \mul \overline{A} = \span \{ \overline {S}e -f\}.
\]
According to the assumption that the sum \eqref{defe} is direct
there are two cases for the relation $A$ not to be regular. In the
first case $e \in \dom S$ and $Se \neq f$ and \eqref{defe1} is
satisfied; this actually means that $\mul A$ is not trivial and
$A$ itself is not an operator. In the second case $e \not\in \dom
S$ and \eqref{defe1} means $e \in \dom \overline{S} \setminus \dom
S$ and $\overline {S}e \ne f$; in this case $A$ is an operator and
$\mul \overline{A}$ is one-dimensional.

In other words, under the assumption that the sum \eqref{defe} is
direct, the relation $A$ is regular if and only if $e \in \dom
\overline{S} \setminus \dom S$ and $\overline{S}e=f$, or
$e\not\in\dom \overline{S}$.
\end{example}

\begin{example}
\label{ut} Let $A$ be an operator in a Hilbert space $\sH$ which
is not closable with $\mul \overline{A}=\span\{\varphi\}$,
$\|\varphi\|=1$, cf. Example \ref{StSz}. Denote the orthogonal
projection onto $\mul \overline{A}$ by $P$. Then
\[
Ah=(I-P)Ah+PAh, \quad h \in \dom A,
\]
and
\[
PA h=(PAh,\varphi)\varphi = (Ah,\varphi)\varphi, \quad h \in \dom
A.
\]
Hence, by definition,
\[
  A_{\rm{reg}}h=Ah-(Ah,\varphi)\varphi, \quad
  A_{\rm{sing}}h=(Ah,\varphi)\varphi, \quad h \in \dom A,
\]
and by Theorem \ref{jor} $A_{\rm{reg}}$ is regular and
$A_{\rm{sing}}$ is singular. According to Proposition \ref{sin} one
has
\[
(A_{\rm{sing}})^-=\cdom A \times \span\{\varphi\}.
\]
\end{example}

\begin{example}
Let $A$ be an operator which is not closable with $\mul
\overline{A}=\span\{\varphi\}$, $\|\varphi\|=1$, cf. Example
\ref{StSz}, and let $B$ be a bounded everywhere defined operator.
Define the operator $T$ as an operator sum by
\[
  T=A+B, \quad \dom T=\dom A.
\]
Then $T^*=A^*+B^*$ and
\[
  \overline{T}=\overline{A}+B, \quad \dom \overline{T}=\dom
\overline{A}.
\]
Moreover, $\mul \overline{T}=\mul \overline{A}=\span \{\varphi\}$,
which implies that the operator $T=A+B$ is not closable. It follows
from \eqref{oper1} and \eqref{oper2} that
\[
  T_{\rm{reg}}h=[A_{\rm{reg}} +(I-P)B]h, \quad
  T_{\rm{sing}}h=A_{\rm{sing}}h+ (Bh,\varphi)\varphi, \quad h \in
  \dom A.
\]
Again from Proposition \ref{sin} one obtains
\[
(T_{\rm{sing}})^-=\cdom A \times \span\{\varphi\}.
\]
In particular, if the operator $A$ is singular, then its
perturbation by the bounded operator $B\neq 0$ is singular if and
only if $\ran B=\span\{\varphi\}$.
\end{example}

\begin{example}
Let $A$ be a bounded operator in a Hilbert space $\sH$ with domain
$\dom A$, let $\sR$ be a not necessarily closed subspace of $\sH$,
and define the linear relation $T$ by
\[
  T=\{\,\{f,Af+\varphi\} :\,f \in \dom A, \, \varphi \in \sR\,\}.
\]
It is not difficult to see that
\[
  \overline{T}=\{\,\{f,\overline{A}f+\varphi\} :\,f \in \dom
\overline{A}, \, \varphi \in {\sR}^-\,\},
\]
since $\overline{A}$ is a bounded operator, cf. \eqref{summ} and
\eqref{summ1}. Now this identity implies that $\mul
\overline{T}={\sR}^-$. Hence
\[
  T_{\rm reg}=\{\,\{f, (I-P)Af\}:\, f \in \dom A\,\},
\]
and
\[
  T_{\rm sing}=\{\,\{f, PAf +\varphi\}:\, f \in \dom A,\, \varphi
  \in \sR\,\}.
\]
It follows from Proposition \ref{sin} that
\[
  (T_{\rm sing})^- = \cdom A \times {\sR}^-.
\]
Note that $T$ is maximally singular if and only if $\sR$ is dense
in $\sH$. If in this example $A$ is the null operator on $\sD=\dom
A$, then $T$ reduces to the Cartesian product
\[
 T=\sD \times \sR,
\]
 which is a singular relation.
\end{example}

\section{An approach via adjoint relations}

There is another approach to the decomposition results in Section
4. It is based on a description of the domains of the adjoints of
$T_{\rm reg}$ and $T_{\rm sing}$. First a general characterization
result will be described.

Let $T$ be a linear relation from a Hilbert space $\sH$ to a Hilbert
space $\sK$. The following lemma gives a description of $\dom T^*$.
Clearly, by inverting the linear relation, a similar result can be
obtained for the description of $\ran T^*$. This latter description
goes back to \cite{Sh} in the case of bounded linear operators (see
also \cite{FW}) and to \cite{Se} for densely defined operators; the
present version can be found in \cite{H}, the proof is presented
here for completeness. The orthogonal decomposition of a closed
linear relation $K$ into an operator part and a multivalued part is
already assumed here (see Proposition~\ref{ca}):
\begin{equation}
\label{odecom}
 K=K_s \hoplus K_\infty,
\end{equation}
where $K_\infty\okr\{0\}\times \mul K$ and $K_s\okr K\hominus
K_\infty\,(=K_{\rm reg})$.

\begin{lemma} {\rm (\cite{H})}
\label{lemma1} Let $T$ be a linear relation from a Hilbert space
$\sH$ to a Hilbert space $\sK$. Then $g\in\dom T^*$ if and only if
there exists a nonnegative number $C_g$, such that
\begin{equation}
\label{est01} |(f',g)_\sK|\le C_g \|f\|_\sH \quad \mbox{for all}
\quad \{f,f'\} \in T.
\end{equation}
In this case the smallest $C_g$ satisfying \eqref{est01} is
$C_g=\|g'\|_\sH$ with $\{g,g'\}\in T^*$ and $g'\in\cdom T$, i.e.,
$C_g=\|(T^*)_s g\|_\sH$.
\end{lemma}

\begin{proof}
First assume that $g\in\dom T^*$. Then $\{g,g'\}\in T^*$ for some
$g'\in\sH$ and by the definition of the adjoint $T^*$ one obtains
for every $\{f,f'\}\in T$:
\[
 |(f',g)_\sK|=|(f,g')_\sH|\le \|f\|_\sH \|g'\|_\sH,
\]
so that one can take $C_g=\|g'\|_\sH$ in \eqref{est01}.

Conversely, assume that $g\in\sK$ satisfies the estimate
\eqref{est01}. Define the linear relation $L_g$ in $\sH\oplus\dC$
by
\[
 L_g:=\left\{\,\{f,(f',g)_\sK\}:\, \{f,f'\}\in T \,\right\}.
\]
Then it follows from \eqref{est01} that $L_g$ is single-valued,
since $\|f\|_\sH=0$ implies $(f',g)_\sK=0$. Hence, $L_g$ is (the
graph of) a single-valued bounded linear functional defined on
$\dom T$. Therefore, it has a continuation $\bar{L}_g$ from $\cdom
T$ into $\dC$ with the same norm ($\|\bar{L}_g\|\le C_g$). By the
Riesz Representation Theorem there exists $g'\in\cdom T$ with
$\|g'\|_\sH=\|\bar{L}_g\|$, such that
\[
 \bar{L}_gf=(f,g')_\sH \quad \mbox{for all } f\in\cdom T.
\]
Therefore $(f',g)_\sK=(f,g')_\sH$ holds for every $\{f,f'\}\in T$,
so that $\{g,g'\}\in T^*$. In particular, $g\in\dom T^*$.

The last statement is clear from the given arguments and the
definition of the orthogonal operator part $(T^*)_s$ of $T^*$, cf.
\eqref{odecom}.
\end{proof}

Let $T$ be a linear relation from a Hilbert space $\sH$ to a Hilbert
space $\sK$ and let the regular and singular parts $T_{\rm reg}$ and
$T_{\rm sing}$ of $T$ be as defined in \eqref{oper1} and
\eqref{oper2}, respectively. In the present approach $\dom (T_{\rm
reg})^*$, $\ker (T_{\rm reg})^*$, $\dom (T_{\rm sing})^*$, and $\ker
(T_{\rm sing})^*$ are determined by means of Lemma~\ref{lemma1}.

\begin{proposition}
\label{adjprop} Let $T$ be a linear relation from a Hilbert space
$\sH$ into a Hilbert space $\sK$. Then
\begin{equation}
\label{newa1}
  \dom (T_{\rm reg})^*=\dom T^* \oplus \mul \overline{T},
  \quad
  \ker (T_{\rm reg})^*=\ker T^* \oplus \mul \overline{T},
\end{equation}
and
\begin{equation}
\label{newa2}
  \dom (T_{\rm sing})^*=\ker (T_{\rm sing})^*=\cdom T^*.
\end{equation}
\end{proposition}

\begin{proof}
It follows from the definition of $T_{\rm reg}$ in \eqref{oper1}
and from Lemma \ref{lemma1} that
\begin{equation}
\label{kar1} g \in \dom (T_{\rm reg})^* \quad \Leftrightarrow
\quad |(f',(I-P)g)| \le C_g \|f\|, \quad \{f,f'\} \in T.
\end{equation}
Clearly, the following equivalence is a consequence of
\eqref{kar1}:
\[
 g \in \dom (T_{\rm reg})^* \cap (\mul \overline{T})^\perp \quad
 \Leftrightarrow \quad g \in \dom T^* \cap (\mul
 \overline{T})^\perp.
\]
Therefore
\[
  \dom (T_{\rm reg})^* \cap \cdom T^*=\cdom T^* \cap \dom T^*,
\]
or, equivalently,
\begin{equation}
\label{kar2}
 \dom (T_{\rm reg})^* \cap \cdom T^*= \dom T^*.
\end{equation}
Furthermore it follows from \eqref{kar1} that
\begin{equation}
\label{kar3}
 \mul \overline{T} \subset \dom (T_{\rm reg})^*.
\end{equation}
Since $\mul \overline{T}\oplus\cdom T^*=\sK$, the first identity in
\eqref{newa1} is obtained from \eqref{kar2} and \eqref{kar3}. In
addition, by the second statement in Lemma~\ref{lemma1} $g\in\ker
T^*(=\ker (T^*)_s)$ if and only if the smallest constant in
\eqref{est01} is $C_g=0$. Therefore, if follows from \eqref{kar1}
that
\begin{equation}
\label{kar2b}
 \ker (T_{\rm reg})^* \cap \cdom T^*= \ker T^* \mbox{ and }
 \mul \overline{T} \subset \ker (T_{\rm reg})^*,
\end{equation}
which proves the second identity in \eqref{newa1}.

Likewise, it follows from the definition of $T_{\rm sing}$ in
\eqref{oper2} and from Lemma \ref{lemma1} that
\begin{equation}
\label{kar} g \in \dom (T_{\rm sing})^* \quad \Leftrightarrow
\quad |(f',Pg)| \le C_g \|f\|, \quad \{f,f'\} \in T.
\end{equation}
Clearly, the following equivalence is a consequence of
\eqref{kar}:
\[
 g \in \dom (T_{\rm sing})^* \cap \mul \overline{T} \quad
 \Leftrightarrow \quad g \in \dom T^* \cap \mul \overline{T}.
\]
Observe that $\dom T^* \cap \mul \overline{T}=\{0\}$, which shows
that
\begin{equation}
\label{kar4}
  \dom (T_{\rm sing})^* \cap \mul \overline{T} =\{0\}.
\end{equation}
Furthermore it follows from \eqref{kar} that
\begin{equation}
\label{kar5} \cdom T^*=(\mul \overline{T})^\perp \subset \ker
(T_{\rm sing})^*\subset \dom (T_{\rm sing})^*
\end{equation}
Since $\mul \overline{T}\oplus\cdom T^*=\sK$, the equalities in
\eqref{newa2} follow from \eqref{kar4} and \eqref{kar5}.
\end{proof}

The results in Proposition~\ref{adjprop} immediately give the main
result in Section~\ref{Sec4}, see Theorem~\ref{jor}. For instance,
according to \eqref{newa1}
\[
 \cdom (T_{\rm reg})^*=\cdom T^* \oplus \mul \overline{T}=\sK,
\]
and thus $T_{\rm reg}$ is regular by Proposition~\ref{pro}. Also the
first equality in \eqref{newa2} shows that $T_{\rm sing}$ is
singular by part (ii) of Proposition~\ref{sin}.

\end{document}